%&amstex
\input amstex
\input amsppt.sty
\magnification=\magstep1
%\textwidth=13.5cm
%\textheight=24cm
%\hoffset=-1cm
%\hsize=30truecc
%\baselineskip=16truept
%\vsize=22.2truecm
\hsize=33truecc
\vsize=22.2truecm
\baselineskip=16truept
\NoBlackBoxes
\nologo
\pageno=1
\topmatter
\TagsOnRight

\def\N{\Bbb N}
\def\Z{\Bbb Z}

\def\l{\left}
\def\r{\right}
\def\b{\bigg}

\def\({\b(}
\def\[{\b[}
\def\){\b)}
\def\]{\b]}

\def\t{\text}
\def\f{\frac}

\def\ord{\roman{ord}}

\def\em{\emptyset}
\def\se {\subseteq}

\def\bi{\binom}
\def\eq{\equiv}

\def\ls{\leqslant}
\def\gs{\geqslant}
\def\al{\alpha}

\def\da{\delta}

\def\Proof{\noindent{\it Proof}}
\def\Remark{\noindent{\it Remark}}

\hbox{Preprint, {\tt arXiv:1504.01608}}
\medskip
\title Natural numbers represented by $\lfloor x^2/a\rfloor+\lfloor y^2/b\rfloor+\lfloor z^2/c\rfloor$\endtitle
\author Zhi-Wei SUN \endauthor
\affil Department of Mathematics, Nanjing University
     \\Nanjing 210093, People's Republic of China
    \\  zwsun\@nju.edu.cn
    \\ {\tt http://math.nju.edu.cn/$\sim$zwsun}
 \endaffil
\abstract Let $a,b,c$ be positive integers. It is known that there are infinitely many positive integers not representated by $ax^2+by^2+cz^2$ with $x,y,z\in\Z$.
In contrast, we conjecture that any natural number is represented by $\lfloor x^2/a\rfloor+\lfloor y^2/b\rfloor +\lfloor z^2/c\rfloor$
with $x,y,z\in\Z$ if $(a,b,c)\not=(1,1,1),(2,2,2)$, and that any natural number is represented by $\lfloor T_x/a\rfloor+\lfloor T_y/b\rfloor +\lfloor T_z/c\rfloor$
with $x,y,z\in\Z$, where $T_x$ denotes the triangular number $x(x+1)/2$. We confirm this general conjecture in some special cases; in particular, we prove that
$$\left\{x^2+y^2+\left\lfloor\frac{z^2}5\right\rfloor:\ x,y,z\in\Z\ \t{and}\ 2\nmid y\right\}=\{1,2,3,\ldots\}$$
and
$$\left\{\left\lfloor\frac{x^2}m\right\rfloor+\left\lfloor\f{y^2}m\right\rfloor+\left\lfloor\frac{z^2}m\right\rfloor:\ x,y,z\in\Z\right\}
=\{0,1,2,\ldots\}\quad \t{for}\ m=5,6,15.$$
We also pose several conjectures for further research; for example, we conjecture that any integer can be written as $x^4-y^3+z^2$, where $x$, $y$ and $z$ are positive integers.
\endabstract
\thanks 2010 {\it Mathematics Subject Classification}.
Primary 11E25; Secondary 11B75, 11D85, 11E20, 11P32.
\newline \indent Supported by the National Natural Science
Foundation (Grant No. 11171140) of China.
\endthanks
\endtopmatter
\document

\heading{1. Introduction}\endheading

Let $\N=\{0,1,2,\ldots\}$ be the set of all natural numbers (nonnegative integers).
A well-known theorem of Lagrange asserts that each $n\in\N$ can be written as the sum of four squares.
It is known that for any $a,b,c\in\Z^+=\{1,2,3,\ldots\}$ there are infinitely many positive integers not represented by $ax^2+by^2+cz^2$
with $x,y,z\in\Z$.

A classical theorem of Gauss and Legendre states that $n\in\N$  can be written as the sum of
three squares if and only if it is not of the form $4^k(8l+7)$ with $k,l\in\N$. Consequently, for each $n\in\N$ there are $x,y,z\in\Z$ such that
$$8n+3=(2x+1)^2+(2y+1)^2+(2z+1)^2,\ \t{i.e.}, \ n=\f{x(x+1)}2+\f{y(y+1)}2+\f{z(z+1)}2.$$
Those $T_x=x(x+1)/2$ with $x\in\Z$ are called {\it triangular numbers}.
For $m=3,4,\ldots$, those {\it $m$-gonal numbers} (or {\it polygonal numbers of order $m$}) are given
by
$$p_m(n):=(m-2)\bi n2+n=\f{(m-2)n^2-(m-4)n}2\ (n=0,1,2,\ldots),$$
and those $p_m(x)$ with $x\in\Z$ are called {\it generalized $m$-gonal numbers}.
Cauchy's polygonal number theorem states that for each $m=5,6,\ldots$ any $n\in\N$
can be written as the sum of $m$ polygonals of order $m$ (see, e.g., [N96, pp.\, 3-35] and [MW, pp.\,54-57].)

For any $k\in\Z$, we clearly have
$$T_k=\f{(2k+1)^2-1}8=\l\lfloor\f{(2k+1)^2}8\r\rfloor.$$
As any natural number can be expressed as the sum of three triangular numbers,
each $n\in\N$ can be written as $\lfloor x^2/8\rfloor+\lfloor y^2/8\rfloor+\lfloor z^2/8\rfloor$ with $x,y,z\in\Z$.
B. Farhi [F13] conjectured that any $n\in\N$ can be expressed the sum of three elements of the set $\{\lfloor x^2/3\rfloor:\ x\in\Z\}$
and showed this for $n\not\eq 2\pmod{24}$. The conjecture was later proved by S. Mezroui, A. Azizi and M. Ziane [MAZ]
in 2014 via the known formula for the number of ways to write $n$ as the sum of three squares. In [F] Farhi provided an elementary proof of the conjecture
and made a further conjecture that for each $a=3,4,5,\ldots$ any $n\in\N$ can be written as
the sum of three elements of the set $\{\lfloor x^2/a\rfloor:\ x\in\Z\}$.
This general conjecture of Farhi has been solved for $a=3,4,7,8,9$ (cf. [HKR]).

Motivated by the above work, we pose the following general conjecture based on our computation.

\proclaim{Conjecture 1.1} Let $a,b,c\in\Z^+$ with $a\ls b\ls c$.

{\rm (i)} If the triple $(a,b,c)$ is neither $(1,1,1)$ nor $(2,2,2)$, then for any $n\in\N$ there are $x,y,z\in\Z$ such that
$$n=\l\lfloor \f{x^2}a\r\rfloor+\l\lfloor \f{y^2}b\r\rfloor+\l\lfloor \f{z^2}c\r\rfloor
=\l\lfloor \f{x^2}a+\f{y^2}b\r\rfloor+\l\lfloor \f{z^2}c\r\rfloor
=\l\lfloor \f{x^2}a\r\rfloor+\l\lfloor \f{y^2}b+\f{z^2}c\r\rfloor.\tag1.1$$

{\rm (ii)} For any $n\in\N$, there are $x,y,z\in\Z$ such that
$$n=\l\lfloor\f{T_x}a\r\rfloor+\l\lfloor \f{T_y}b\r\rfloor+\l\lfloor \f{T_z}c\r\rfloor
=\l\lfloor\f{T_x}a+\f{T_y}b\r\rfloor+\l\lfloor \f{T_z}c\r\rfloor
=\l\lfloor\f{T_x}a\r\rfloor+\l\lfloor \f{T_y}b+\f{T_z}c\r\rfloor.\tag1.2$$
Moreover, if the triple $(a,b,c)$ is not among
$$(1,1,1),\ (1,1,3),\ (1,1,7),\ (1,3,3),\ (1,7,7),\ (3,3,3),$$
then for any $n\in\N$ there are $x,y,z\in\Z$ such that
$$\aligned n=&\l\lfloor\f{x(x+1)}a\r\rfloor+\l\lfloor \f{y(y+1)}b\r\rfloor+\l\lfloor \f{z(z+1)}c\r\rfloor
\\=&\l\lfloor\f{x(x+1)}a+\f{y(y+1)}b\r\rfloor+\l\lfloor \f{z(z+1)}c\r\rfloor
\\=&\l\lfloor\f{x(x+1)}a\r\rfloor+\l\lfloor \f{y(y+1)}b+\f{z(z+1)}c\r\rfloor.
\endaligned\tag1.3$$
\endproclaim

In this paper we establish some results in the direction of Conjecture 1.1.

\proclaim{Theorem 1.1} {\rm (i)} For each $m=4,6$, any $n\in\N$ can be written as $x^2+(2y)^2+\lfloor z^2/m\rfloor$ with $x,y,z\in\Z$.
Also, any $n\in\Z^+$ can be expressed as $x^2+y^2+\lfloor z^2/5\rfloor$ with $x,y,z\in\Z$ and $2\nmid y$.

{\rm (ii)} For any $\da\in\{0,1\}$, any $n\in\Z^+$ can be expressed as $x^2+y^2+\lfloor z^2/8\rfloor$ with $x,y,z\in\Z$ and $y\eq \da\pmod 2$.

{\rm (iii)} For each $m=2,3,9,21$, any $n\in\N$ can be written as $x^2+y^2+\lfloor z^2/m\rfloor$ with $x,y,z\in\Z$.
Also, for each $m=3,4,6$ we have
$$\l\{x^2+y^2+\l\lfloor\f{z(z+1)}m\r\rfloor:\ x,y,z\in\Z\r\}=\N.\tag1.4$$

{\rm (iv)} For each $m=5,6,15$, we have
$$\l\{x^2+\l\lfloor\f{y^2}{m}\r\rfloor+\l\lfloor\f{z^2}{m}\r\rfloor:\, x,y,z\in\Z\r\}
=\l\{\l\lfloor\f{x^2}{m}\r\rfloor+\l\lfloor\f{y^2}{m}\r\rfloor+\l\lfloor\f{z^2}{m}\r\rfloor:\, x,y,z\in\Z\r\}=\N.\tag1.5$$

{\rm (v)} We have
$$\l\{T_x+T_y+\l\lfloor\f{T_z}3\r\rfloor:\ x,y,z\in\Z\r\}=\l\{\l\lfloor\f{T_x}3\r\rfloor+\l\lfloor\f{T_y}3\r\rfloor+\l\lfloor\f{T_z}3\r\rfloor:\ x,y,z\in\Z\r\}=\N$$
and
$$\l\{x(x+1)+y(y+1)+\l\lfloor\f{z(z+1)}4\r\rfloor:\ x,y,z\in\Z\r\}=\N.\tag1.6$$
\endproclaim
\Remark\ 1.1. As $x^2=(3x)^2/9$, Theorem 1.1(iii) with $m=9$ implies that
any $n\in\N$ can be written as $\lfloor x^2/9\rfloor+\lfloor y^2/9\rfloor+\lfloor z^2/9\rfloor$ with $x,y,z\in\Z$.
Theorem 1.1(iv) confirms Farhi's conjecture for $a=5,6,15$. The author [S15a, Remark 1.8] conjectured that for any $n\in\N$
we can write $20n+9$ as $5x^2+5y^2+(2z+1)^2$ with $x,y,z\in\Z$; it is easy to see that (1.4) for $m=5$ follows from this conjecture.
As $\{2T_x+2T_y+T_z:\ x,y,z\in\Z\}=\N$ by Liouville's result, any $n\in\N$ can be written as $T_x+T_y+T_z/2$ with $x,y,z\in\Z$.

\medskip

As a supplement to parts (i)-(iii) of Theorem 1.1, we pose the following conjecture.

\proclaim{Conjecture 1.2} {\rm (i)} Let $n\in\Z^+$. Then, for any integer $m>6$ and $\da\in\{0,1\}$, we have
$n=x^2+y^2+\lfloor z^2/m\rfloor$ for some $x,y,z\in\Z$ with $y\eq \da\pmod 2$.

{\rm (ii)} For any integer $m>2$, we have
$$\l\{x^2+(2y)^2+\l\lfloor\f{z(z+1)}m\r\rfloor:\ x,y,z\in\Z\r\}=\N.$$
For each $m=4,5,\ldots$, any positive integer $n$ can be represented by $x^2+y^2+\lfloor z(z+1)/m\rfloor$ with $x,y,z\in\Z$ and $2\nmid y$.
\endproclaim
\Remark\ 1.2. It is known that $\{x^2+(2y)^2+T_z:\ x,y,z\in\Z\}=\{x^2+(2y)^2+2T_z:\ x,y,z\in\Z\}=\N$ (cf. [S07, Section 4]).
\medskip

For any $a\in\Z^+$, clearly
$$\l\{\l\lfloor \f{x^2}a\r\rfloor:\ x\in\Z\r\}\supseteq\l\{\l\lfloor\f{(ax)^2}a\r\rfloor=ax^2:\ x\in\Z\r\}.$$

\proclaim{Theorem 1.2} {\rm (i)} For each $m=2,3,4,5$ we have
$$\l\{x^2+2y^2+\l\lfloor\f{z^2}m\r\rfloor:\ x,y,z\in\Z\r\}=\N.$$

{\rm (ii)} For each $m=3,4,6,8$, we have
$$\l\{x^2+3y^2+\l\lfloor\f{z^2}m\r\rfloor:\ x,y,z\in\Z\r\}=\N.$$

{\rm (iii)} We have
$$\l\{x^2+5y^2+\l\lfloor\f{z^2}8\r\rfloor:\ x,y,z\in\Z\r\}=\l\{x^2+6y^2+\l\lfloor\f{z^2}4\r\rfloor:\ x,y,z\in\Z\r\}=\N.$$

{\rm (iv)} We have
$$\l\{2x^2+2y^2+\l\lfloor\f{z^2}8\r\rfloor:\ x,y,z\in\Z\r\}=\l\{2x^2+3y^2+\l\lfloor\f{z^2}3\r\rfloor:\ x,y,z\in\Z\r\}=\N$$
and
$$\l\{2x^2+\l\lfloor\f{y^2}2\r\rfloor+\l\lfloor\f{z^2}3\r\rfloor:\ x,y,z\in\Z\r\}=\N.$$
\endproclaim

Our following conjecture involving the ceiling function is quite similar to Conjecture 1.1.

\proclaim{Conjecture 1.3} Let $a,b,c\in\Z^+$ with $a\ls b\ls c$.

{\rm (i)} If the triple $(a,b,c)$ is not among $(1,1,1),(1,1,2),(1,1,5)$, then for any $n\in\N$ there are $x,y,z\in\Z$ such that
$$n=\l\lceil \f{x^2}a\r\rceil+\l\lceil \f{y^2}b\r\rceil+\l\lceil \f{z^2}c\r\rceil.$$

{\rm (ii)} We have
$$\l\{\l\lceil\f{T_x}a\r\rceil+\l\lceil \f{T_y}b\r\rceil+\l\lceil\f{T_z}c\r\rceil:\ x,y,z\in\Z\r\}=\N.$$
Moreover, if the triple $(a,b,c)$ is neither $(1,1,1)$ nor $(1,1,3)$,
then for any $n\in\N$ there are $x,y,z\in\Z$ such that
$$n=\l\lceil\f{x(x+1)}a\r\rceil+\l\lceil \f{y(y+1)}b\r\rceil+\l\lceil \f{z(z+1)}c\r\rceil.$$
\endproclaim

We are also able to deduce some results similar to Theorems 1.1-1.2 in the direction of Conjecture 1.3.
Here we just collect few results of this type.

\proclaim{Theorem 1.3} {\rm (i)} For each $m=2,3,4,5,6,15$, we have
$$\l\{\l\lceil\f{x^2}m\r\rceil+\l\lceil\f{y^2}m\r\rceil+\l\lceil\f{z^2}{m}\r\rceil:\ x,y,z\in\Z\r\}=\N.\tag1.7$$

{\rm (ii)} We have
$$\l\{x^2+3y^2+\l\lceil\f{z^2}{2}\r\rceil:\ x,y,z\in\Z\r\}=\l\{x^2+3y^2+\l\lceil\f{z^2}{10}\r\rceil:\ x,y,z\in\Z\r\}=\N.\tag1.8$$

{\rm (iii)} For any $n\in\N$, there are $x,y,z\in\Z$ such that
$$n=x(x+1)+\f{y(y+1)}3+\l\lceil\f{z(z+1)}3\r\rceil.\tag1.9$$
Also, any $n\in\N$ can be written as $x(3x+1)+y(3y+1)+\lceil z(z+1)/3\rceil$ with $x,y,z\in\Z$, and hence
$$\l\{\l\lceil \f{x(x+1)}3\r\rceil+\l\lceil\f{y(y+1)}3\r\rceil+\l\lceil\f{z(z+1)}3\r\rceil:\ x,y,z\in\Z\r\}=\N.\tag1.10$$
\endproclaim
\Remark\ 1.3. In contrast with (1.8), we note that $20142$ is the first natural number not represented by $x^2+3y^2+\lfloor z^2/10\rfloor$ with $x,y,z\in\Z$.
\medskip

Now we state another theorem.

\proclaim{Theorem 1.4} {\rm (i)} For any integer $a>1$, we have
$$\l\{\l\lfloor\f{x^2+y^2+z^2}a\r\rfloor:\ x,y,z\in\Z\r\}=\N$$
and
$$\l\{\l\lfloor\f{x(x+1)+y(y+1)+z(z+1)}a\r\rfloor:\ x,y,z\in\Z\r\}=\N.$$

{\rm (ii)} Let $a\in\Z^+$. If $a$ is odd, then any $n\in\N$ can be written as
$x^2+y^2+z^2+\lfloor\f a2(x+y+z)\rfloor$ with $x,y,z\in\Z$. If $3\nmid a $, then any $n\in\N$ can be written as
$x^2+y^2+z^2+\lfloor\f a3(x+y+z)\rfloor$ with $x,y,z\in\Z$.

{\rm (iii)} For any $n\in\N$, there are $x,y,z\in\Z$ such that
$$n=\f{p_8(x)}2+\l\lceil\f{p_8(y)}2\r\rceil+\l\lceil\f{p_8(z)}2\r\rceil.$$
Hence
$$\{s(x)+s(y)+s(z):\ x,y,z\in\Z\}=\N,\tag1.11$$
where
$$s(x):=\l\lceil\f{p_8(-x)}2\r\rceil=x+\l\lceil 1.5x^2\r\rceil.$$
\endproclaim
\Remark\ 1.4. For $m=19,20$, we have $111\not=x^2+y^2+z^2+\lfloor (x+y+z)/m\rfloor$ for any $x,y,z\in\Z$.
\medskip

The generalized octagonal numbers $p_8(x)=x(3x-2)\ (x\in\Z)$ have some properties similar to certain properties of squares.
For example, recently the author [S16] showed that any $n\in\N$ can be written as the sum of four generalized octagonal numbers;
this result is quite similar to Lagrange's theorem on sums of four squares. Note that
$$\l\lfloor \f{p_8(x)}{2m}\r\rfloor=\l\lfloor \f{4p_8(x)+1}{4m}\r\rfloor=\l\lfloor\f{p_8(1-2x)}{4m}\r\rfloor
\ \t{and}\ \l\lfloor \f{p_8(x)}m\r\rfloor=\l\lfloor \f{(3x-1)^2}{3m}\r\rfloor\tag1.12$$
 for any $m\in\Z^+$ and $x\in\Z$.

\proclaim{Theorem 1.5} {\rm (i)} $n\in\N$ can be written as $p_8(x)+p_8(y)+2p_8(z)$ with $x,y,z\in\Z$ if and only if $n$
 does not belong to the set
 $$\l\{4^{k+2}q-\f23(4^k+2):\ k\in\N\ \t{and}\ q\in\Z^+\r\}.$$
Also, each nonnegative even number  can be represented by $p_8(x)+2p_8(y)+4p_8(z)$ with $x,y,z\in\Z$. Consequently,
$$\l\{p_8(x)+\l\lfloor\f{p_8(y)}2\r\rfloor+\l\lfloor\f{p_8(z)}2\r\rfloor:\ x,y,z\in\Z\r\}=\N\tag1.13$$
and $$\l\{\l\lfloor\f{x^2}3\r\rfloor+\l\lfloor\f{y^2}6\r\rfloor+\l\lfloor\f{z^2}6\r\rfloor:\ x,y,z\in\Z\r\}=\N.\tag1.14$$

{\rm (ii)} We have
$$\l\{p_8(x)+p_8(y)+\l\lfloor\f{p_8(z)}2\r\rfloor:\ x,y,z\in\Z\r\}=\N,\tag1.15$$
hence
$$\l\{p_8(x)+p_8(y)+\l\lfloor\f{p_8(z)}8\r\rfloor:\ x,y,z\in\Z\r\}=\N\tag1.16$$
and
$$\l\{\l\lfloor\f{x^2}3\r\rfloor+\l\lfloor\f{y^2}3\r\rfloor+\l\lfloor\f{z^2}6\r\rfloor:\ x,y,z\in\Z\r\}=\N.\tag1.17$$

{\rm (iii)} For $n\in\N$ there are $x,y,z\in\Z$ such that
$$n=p_8(x)+p_8(y)+\f{p_8(z)}4.\tag1.18$$

{\rm (iv)}  We have
$$\l\{p_8(x)+p_8(y)+\l\lfloor\f{p_8(z)}5\r\rfloor:\ x,y,z\in\Z\r\}=\N\tag1.19$$
and hence
$$\l\{\l\lfloor\f{x^2}3\r\rfloor+\l\lfloor\f{y^2}3\r\rfloor+\l\lfloor\f{z^2}{15}\r\rfloor:\ x,y,z\in\Z\r\}=\N.\tag1.20$$
\endproclaim

We are going to prove Theorems 1.1-1.2  in the next section, and show Theorems 1.3-1.4 in Section 3.
Section 4 is devoted to our proof of Theorem 1.5. We pose some further conjectures in Section 5.

\heading{2. Proofs of Theorems 1.1-1.2}\endheading

\proclaim{Lemma 2.1} Suppose that $n\in\Z^+$ is not a power of two. Then there are $x,y,z\in\Z$ with $|x|<n$, $|y|<n$ and $|z|<n$
such that $x^2+y^2+z^2=n^2$.
\endproclaim
\Proof. In 1907 Hurwitz (cf. [D99, p.\,271]) showed that
$$\aligned&|\{(x,y,z)\in\Z^3:\ x^2+y^2+z^2=n^2\}|
\\=&6\prod_{p>2}\l(\f{p^{\ord_p(n)+1}-1}{p-1}+(-1)^{(p+1)/2}\f{p^{\ord_p(m)}-1}{p-1}\r),
\endaligned\tag2.1$$
where $\ord_p(n)$ is the order of $n$ at the prime $p$. Note that
$$(\pm n)^2+0^2+0^2=0^2+(\pm n)^2+0^2=0^2+0^2+(\pm n)^2.$$
As $n$ has an odd prime $p$, by (2.1) we have
$$|\{(x,y,z)\in\Z^3:\ x^2+y^2+z^2=n^2\}|\gs 6\f{p^{\ord_p(n)+1}-p^{\ord_p(n)}}{p-1}\gs 6p>8$$
and hence there are $x,y,z\in\Z$ with $x^2,y^2,z^2\not=n^2$ such that $x^2+y^2+z^2=n^2$.
This concludes the proof. \qed

\proclaim{Lemma 2.2} {\rm (i)} Let $u$ and $v$ be integers with $u^2+v^2$ a positive multiple of $5$.
Then $u^2+v^2=x^2+y^2$ for some $x,y\in\Z$ with $5\nmid xy$.

{\rm (ii)} For any $n\in\N$ with $n\eq \pm6\pmod{20}$, we can write $n$ as $5x^2+5y^2+z^2$ with $x,y,z\in\Z$ and $2\nmid z$.
\endproclaim
\Remark\ 2.1. Parts (i) and (ii) of Lemma 2.2 are Lemmas 2.1 and 2.2 of [S15b].

\proclaim{Lemma 2.3} Let $n>1$ be an integer with $n\eq1,9\pmod{20}$ or $n\eq11,19\pmod{40}$.
Then we can write $n$ as $5x^2+5y^2+z^2$ with $x,y,z\in\Z$ such that $x\not\eq y\pmod 2$ if $n\eq1,9\pmod{20}$,
and $2\nmid y$ if $n\eq11,19\pmod{40}$.
\endproclaim
\Proof. As $n\eq1\pmod 4$ or $n\eq3\pmod8$, by the Gauss-Legendre theorem $n$ is the sum of three squares.
As $n$ is not a power of two, in view of Lemma 2.1 we can always write $n$ as $w^2+u^2+v^2$ with $u,v,w\in\Z$ and $w^2,u^2,v^2\not=n$.
Without loss of generality, we assume that $2\nmid w$ and $u\eq v\pmod 2$. Clearly, $u\eq v\eq 0\pmod2$ if $n\eq1\pmod4$.
If $w^2\eq-n\pmod 5$, then $u^2+v^2\eq2n\pmod 5$ and hence $u^2\eq v^2\eq n\pmod 5$.
If $w^2\eq n\pmod 5$, then $u^2+v^2$ is a positive multiple of $5$ and hence by Lemma 2.2 we can write it as $s^2+t^2$, where $s$ and $t$
are integers with $s^2\eq -n\pmod5$ and $t^2\eq n\pmod 5$.
When $n\eq1\pmod 4$, we have $s^2+t^2=u^2+v^2\eq 0\pmod 4$, we have $s\eq t\eq 0\pmod 2$.
If $5\mid w$, then one of $u^2$ and $v^2$ is divisible by $5$ and the other is congruent to $n$ modulo $5$.

By the above, there always exist $x,y,z\in\Z$ with $z^2\eq n\pmod5$ such that $n=x^2+y^2+z^2$ and that
$2\mid z$ if $n\eq1\pmod 4$.
Note that $x^2\eq-y^2\eq(\pm2y)^2\pmod 5$. Without loss of generality, we assume that $x\eq2y\pmod 5$ and hence $2x\eq-y\pmod5$.
Set $\bar x=(x-2y)/5$ and $\bar y=(2x+y)/5$. Then
$$n=x^2+y^2+z^2=5\bar x^2+5\bar y^2+z^2.$$
If $n\eq1\pmod 4$, then $2\mid z$ and hence $\bar x\not\eq\bar y\pmod2$.
If $n\eq3\pmod 8$, then $z^2\not\eq n\pmod 4$ and hence $\bar x$ or $\bar y$ is odd. This concludes the proof. \qed
\medskip

\Remark\ 2.1. Without using Lemma 2.1 and Lemma 2.2(i), the author [S15a, Theorem 1.7(iv)] showed by a different method that
for any integer $n>1$ with $n\eq1,9\pmod{20}$ we can write $n=5x^2+5y^2+(2z)^2$ with $x,y,z\in\Z$ if $n$ is not a square.

\medskip

For convenience, we define
$$E(f(x,y,z)):=\{n\in\N:\ n\not=f(x,y,z)\ \ \t{for any}\ x,y,z\in\Z\}$$
for any function $f:\Z^3\to\N$.

\medskip
\noindent{\it Proof of Theorem 1.1}. Let $n$ be a fixed nonnegative integer.

(i) By Dickson [D39, pp.\,112-113],
$$E(4x^2+16y^2+z^2)=\bigcup_{k\in\N}\{4k+2,4k+3,16k+12\}\cup\{4^k(8l+7):\ k,l\in\N\}.$$
So, there are $x,y,z\in\Z$ such that $4n+1=4x^2+16y^2+z^2$ and hence $n=x^2+(2y)^2+\lfloor z^2/4\rfloor$.

For $r\in\{1,4\}$, if $6n+r=6x^2+24y^2+z^2$ with $x,y,z\in\Z$, then $z^2\eq r\pmod 6$ and $n=x^2+(2y)^2+\lfloor z^2/6\rfloor$.
By Dickson [D39, pp.\,112-113],
$$E(6x^2+24y^2+z^2)=\bigcup_{k\in\N}\{8k+3,8k+5,32k+12\}\cup\{9^k(3l+2):\ k,l\in\N\}.$$
If both $6n+1$ and $6n+4$ belong to this set, then one of them has the form $32k+12$ and hence we get a contradiction
since $32k+12\pm3\not\eq3,5\pmod 8$.

If $n\eq0,1\pmod 4$, then $5n+1\eq1,6\pmod {20}$ and hence by Lemmas 2.2 and 2.3 we have $5n+1=5x^2+5y^2+z^2$ with $x,y,z\in\Z$ and $x\not\eq y\pmod2$,
thus $x$ or $y$ is odd and $n=x^2+y^2+\lfloor z^2/5\rfloor$.
 By Dickson [D39, pp.\,112-113],
$$E(x^2+y^2+5z^2)=\{4^k(8l+3):\ k,l\in\N\}.$$
 If $n\eq2\pmod 4$ or $n\eq7\pmod 8$, then
  there are $x,y,z\in\Z$ such that $n=x^2+y^2+5z^2=x^2+y^2+\lfloor (5z)^2/5\rfloor$
 and one of $x$ and $y$ is odd since $5z^2\eq z^2\not\eq n\pmod 4$.
If $n\eq3\pmod 8$, then $5n+4\eq 19\pmod {40}$ and hence by Lemma 2.3 there are $x,y,z\in\Z$ with $2\nmid y$ such that
 $5n+4=5(x^2+y^2)+z^2$ and hence $n=x^2+y^2+\lfloor z^2/5\rfloor$ with $y$ odd.

(ii) By [D39, pp.\,112-113], there are $x,y,z\in\Z$ such that $8n+1=8x^2+32y^2+z^2$ and hence $n=x^2+(2y)^2+\lfloor z^2/8\rfloor$.

Suppose that $n\in\Z^+$. As conjectured by Sun [S07] and proved by Oh and Sun [OS], there are $x,y,z\in\Z$ with $y$ odd such that
$n=x^2+y^2+T_z$ and hence $n=x^2+y^2+\lfloor (2z+1)^2/8\rfloor$.

(iii) If $2n\eq 6\pmod 8$, then $2n\not\in\{4^k(8l+7):\ k,l\in\N\}$.
If $2n\not\eq6\pmod 8$, then $2n+1\not\in\{4^k(8l+7):\ k,l\in\N\}$.
So, for some $\da\in\{0,1\}$, we have $2n+\da\not\in\{4^k(8l+7):\ k,l\in\N\}$
and hence (by the Gauss-Legendre theorem) $2n+\da=x^2+y^2+z^2$ for some $x,y,z\in\Z$ with $z\eq\da\pmod 2$.
Note that $x\eq y\pmod 2$ and
$$2n+\da=2\l(\f{x+y}2\r)^2+2\l(\f{x-y}2\r)^2+z^2.$$
Therefore,
$$n=\l(\f{x+y}2\r)^2+\l(\f{x-y}2\r)^2+\f{z^2-\da}2=\l(\f{x+y}2\r)^2+\l(\f{x-y}2\r)^2+\l\lfloor\f{z^2}2\r\rfloor.$$

By Dickson [D39, pp.\,112-113],
$$E(3x^2+3y^2+z^2)=\{9^k(3l+2):\ k,l\in\N\}.$$
So, there are $x,y,z\in\Z$ such that $3n+1=3(x^2+y^2)+z^2$ and hence $n=x^2+y^2+\lfloor z^2/3\rfloor$.

Clearly $9n+1\eq 9n+7\pmod 2$ but $9n+1\not\eq 9n+7\pmod 4$. So, for some $r\in\{1,7\}$, we have $9n+r\not\in\{4^k(8l+7):\ k,l\in\N\}$
and hence (by the Gauss-Legendre theorem) there are $x,y,z\in\Z$ such that $9n+r=(3x)^2+(3y)^2+z^2$ and therefore
$n=x^2+y^2+\lfloor z^2/9\rfloor$.

By Dickson [D39, pp.\,112-113],
$$E(21x^2+21y^2+z^2)=\bigcup_{k,l\in\N}\{4^k(8l+7),9^k(3l+2),49^k(7l+3),49^k(7l+5),49^k(7l+6)\}.$$
For each $r=1,4,16$, if $21n+r$ belongs to the above set then it has the form $4^k(8l+7)$ with $k,l\in\N$.
If $$\{21n+1,21n+4,21n+16\}\se\{4^k(8l+7):\ k,l\in\N\},$$ then $21n+4$ and $21n+16$ are even since
$21n+4\not\eq 21n+16\pmod 8$, hence $21n+1\eq7\pmod 8$ and $21n+4\eq 2\pmod 8$ which leads a contradiction.
So, for some $r\in\{1,4,16\}$ and $x,y,z\in\Z$ we have $21n+r=21(x^2+y^2)+z^2$ and hence $n=x^2+y^2+\lfloor z^2/21\rfloor$.

By Dickson [D39, pp.\,112-113],
$$E(12x^2+12y^2+z^2)=\bigcup_{k\in\N}\{(4k+2,\,4k+3\}\cup\{9^k(3l+2):\ k,l\in\N\}.$$
So, for some $x,y,z\in\Z$ we have $12n+1=12(x^2+y^2)+(2z+1)^2$ and hence
$$n=x^2+y^2+\f{z(z+1)}3=x^2+y^2+\l\lfloor\f{z(z+1)}3\r\rfloor.$$
This proves (1.4) for $m=3$.

 By Jones and Pall [JP], there are $x,y,z\in\Z$ such that $16n+1=16x^2+16y^2+(2z+1)^2$ and hence
$$n=x^2+y^2+\f{(2z+1)^2-1}{16}=x^2+y^2+\l\lfloor\f{z(z+1)}4\r\rfloor.$$
This proves (1.4) for $m=4$.

By [S15a, Theorem 1.7(ii)], $n$ can be written as $x^2+y^2+p_5(z)$ with $x,y,z\in\Z$. Note that
$$p_5(z)=\f{z(3z-1)}2=\f{3z(3z-1)}{6}.$$
So (1.4) holds for $m=6$.

(iv) Now we prove (1.5) for $m=5$. By Dickson [D39, pp.\,112-113],
$$E(5x^2+y^2+z^2)=\{4^k(8l+3):\ k,l\in\N\}.$$
As $5n+2\not\eq 5n+4\pmod 4$, for a suitable choice of $r\in\{2,4\}$ we can write $5n+r$ as $5x^2+y^2+z^2$ with $x,y,z\in\Z$.
If $r=2$, then $y^2\eq z^2\eq1\pmod 5$ and hence
$$n=x^2+\f{y^2-1}5+\f{z^2-1}5=x^2+\l\lfloor\f{y^2}5\r\rfloor+\l\lfloor\f{z^2}5\r\rfloor.$$
If $r=4$, then we may assume that $y^2\eq0\pmod5$ and $z^2\eq 4\pmod 5$, hence
$$n=x^2+\f{y^2}5+\f{z^2-4}5=x^2+\l\lfloor\f{y^2}5\r\rfloor+\l\lfloor\f{z^2}5\r\rfloor.$$

If $\{5n+5,5n+6,5n+9\}\se E:=\{4^k(8l+7):\ k,l\in\N\}$, then we must have $5n+6\eq 7\pmod 8$ and hence $5n+9\eq2\pmod 8$ which leads a contradiction.
Thus, by the Gauss-Legendre theorem, for some $r\in\{0,1,4\}$ the number $5n+5+r$ is the sum of three squares. If $5(n+1)+r=m^2$
for some $m\in\Z^+$ which is not a power of two, then by Lemma 2.1 we have $5(n+1)+r=x^2+y^2+z^2$ for some $x,y,z\in\Z$ with $x^2,y^2,z^2\not=5(n+1)+r$.
If $5(n+1)+r=(2^k)^2$ for some $k\in\Z^+$, then $r\in\{1,4\}$, $5(n+1)+(5-r)=4^k+5-2r\eq 5-2r\eq\pm3\pmod 8$
and hence $5(n+1)+(5-r)\not\in E$. So, for a suitable choice of $r\in\{0,1,4\}$,
we can write $5(n+1)+r=x^2+y^2+z^2$ with $x,y,z\in\Z$ and $x^2,y^2,z^2\not=5(n+1)+r$. Clearly, one of $x^2,y^2,z^2$, say $z^2$, is congruent to $r$ modulo $5$.
Then $x^2+y^2$ is a positive multiple of $5$. By Lemma 2.2, $x^2+y^2=\bar x^2+\bar y^2$ for some $\bar x,\bar y\in\Z$ with $5\nmid \bar x\bar y$.
Without loss of generality we may assume that $\bar x^2\eq1\pmod5$ and $\bar y^2\eq4\pmod5$. Therefore,
$$n=\f{\bar x^2-1}5+\f{\bar y^2-4}5+\f{z^2-r}5=\l\lfloor\f{\bar x^2}5\r\rfloor+\l\lfloor\f{\bar y^2}5\r\rfloor+\l\lfloor\f{z^2}5\r\rfloor.$$

Now we show (1.5) for $m=6$. By Dickson [D39, pp.\,112-113],
$$E(6x^2+y^2+z^2)=\{9^k(9l+3):\ k,l\in\N\}.$$
So, there are $x,y,z\in\Z$ such that $6n+4=6x^2+y^2+z^2$. Clearly, exactly one of $y$ and $z$, say $y$, is divisible by $3$.
Note that $y$ and $z$ have the same parity. If $y\eq z\eq0\pmod 2$, then $y^2\eq 0\pmod 6$ and $z^2\eq 4\pmod 6$.
If $y\eq z\eq1\pmod2$, then $y^2\eq3\pmod6$ and $z^2\eq1\pmod 6$. Anyway, we have
$$n=x^2+\f{y^2+z^2-4}6=x^2+\l\lfloor\f{y^2}6\r\rfloor+\l\lfloor\f{z^2}6\r\rfloor.$$

Assume that $n$ is even. Then $6n+9\eq1\pmod 4$ and hence by the Gauss-Legendre theorem and [S16, Lemma 2.2] we can write $6n+9=x^2+y^2+z^2$ with $x,y,z\in\Z$ and $3\nmid xyz$.
Clearly, exactly one of $x,y,z$, say $x$, is odd. Thus $x^2\eq1\pmod 6$ and $y^2\eq z^2\eq 4\pmod6$. Therefore
$$n=\f{x^2-1}6+\f{y^2-4}6+\f{z^2-4}6=\l\lfloor\f{x^2}6\r\rfloor+\l\lfloor\f{y^2}6\r\rfloor+\l\lfloor\f{z^2}6\r\rfloor.$$

Now suppose that $n$ is odd. Then $3n+4\eq1\pmod 6$, and hence by [S16, Lemma 4.3(ii)] we can write $3n+4=x^2+y^2+2z^2$ with $x,y,z\in\Z$ and $3\nmid xyz$.
Without loss of generality, we may assume that $x\eq y\pmod 3$ (otherwise we may use $-y$ to replace $y$).
Clearly, $x\not\eq y\pmod2$. Thus $6n+8=(x+y)^2+(x-y)^2+(2z)^2$ with $(x+y)^2\eq1\pmod6$, $(x-y)^2\eq 3\pmod6$ and $(2z)^2\eq 4\pmod6$.
Therefore
$$n=\f{(x+y)^2-1}6+\f{(x-y)^2-3}6+\f{(2z)^2-4}6=\l\lfloor\f{(x+y)^2}6\r\rfloor+\l\lfloor\f{(x-y)^2}6\r\rfloor+\l\lfloor\f{(2z)^2}6\r\rfloor.$$

Now we prove (1.5) for $m=15$. By Dickson [D39, pp.\,112-113],
$$E(3x^2+y^2+z^2)=\{9^k(9l+6):\ k,l\in\N\}.$$
So, there are $x,y,z\in\Z$ such that $3n+1=3x^2+y^2+z^2$ and hence
$$15n+5=15x^2+(2^2+1^2)(y^2+z^2)=15x^2+(2y-z)^2+(y+2z)^2.$$
As $(2y-z)^2+(y+2z)^2=5(y^2+z^2)$ is a positive multiple of $5$, by Lemma 2.2 there are $u,v\in\Z$ with $5\nmid uv$
such that $(2y-z)^2+(y+2z)^2=u^2+v^2$. Without loss of generality, we assume that $u^2\eq1\pmod 5$ and $v^2\eq4\pmod 5$.
Then $15n+5=15x^2+u^2+v^2$ with $u^2\eq1\pmod{15}$ and $v^2\eq4\pmod{15}$. Therefore
$$n=x^2+\f{u^2-1}{15}+\f{v^2-1}{15}=x^2+\l\lfloor\f{u^2}{15}\r\rfloor+\l\lfloor\f{v^2}{15}\r\rfloor.$$

If $\{15n+6,15n+9,15n+15\}\se E:=\{4^k(8l+7):\ k,l\in\N\}$, then we must have $15n+6\eq 7\pmod 8$ and hence $15n+9\eq2\pmod 8$ which leads a contradiction.
Thus, by the Gauss-Legendre theorem, for some $r\in\{1,4,10\}$ the number $15n+5+r$ is the sum of three squares.
In view of [S16, Lemma 2.2], we can write $15n+5+r=x^2+y^2+z^2$ with $x,y,z\in\Z$ and $3\nmid xyz$. It is easy to see that
one of $x^2,y^2,z^2$, say $z^2$,  is congruent to $r$ modulo $5$. Then $x^2+y^2$ is a positive multiple of $5$, and hence
by Lemma 2.2 we can write $x^2+y^2=\bar x^2+\bar y^2$ with $\bar x,\bar y\in\Z$ and $5\nmid \bar x\bar y$.
Without loss of generality, we may assume that $\bar x^2\eq1\pmod5$ and $\bar y^2\eq4\pmod5$. Then
$\bar x^2\eq1\pmod{15}$, $\bar y^2\eq4\pmod{15}$ and $z^2\eq r\pmod{15}$. Therefore
$$n=\f{\bar x^2-1}{15}+\f{\bar y^2-4}{15}+\f{z^2-r}{15}=\l\lfloor\f{\bar x^2}{15}\r\rfloor+\l\lfloor\f{\bar y^2}{15}\r\rfloor+\l\lfloor\f{z^2}{15}\r\rfloor.$$

(v) Clearly,
$$\l\{\l\lfloor\f{T_x}3\r\rfloor:\ x\in\Z\r\}\supseteq\l\{p_5(x)=\f{T_{3x-1}}3:\ x\in\Z\r\}.$$
By [S15a, Theorem 1.14], $\{T_x+T_y+p_5(z):\ x,y,z\in\Z\}=\N$.
It is also known that $\{p_5(x)+p_5(y)+p_5(z):\ x,y,z\in\Z\}=\N$ (cf. Guy [Gu] and [S15a]).

Now it remains to prove (1.6). Clearly, for some $r\in\{1,2\}$, $2n+r$ is not a triangular number.
Hence, by [S07, Theorem 1(iii)] there are $x,y,z\in\Z$ with $x\not\eq y\pmod2$ such that
$2n+r=x^2+y^2+T_z$. Thus $4n+2r=(x+y)^2+(x-y)^2+z(z+1)$ with $x\pm y$ odd and $z(z+1)\eq2(r-1)\pmod4$.
Write $x+y=2u+1$ and $x-y=2v+1$ with $u,v\in\Z$. Then
$$\align n=&\f{(2u+1)^2-1}4+\f{(2v+1)^2-1}4+\f{z(z+1)-2(r-1)}4
\\=&u(u+1)+v(v+1)+\l\lfloor\f{z(z+1)}4\r\rfloor.
\endalign$$

In view the above, we have completed the proof of Theorem 1.1. \qed

\medskip
\noindent{\it Proof of Theorem 1.2}. Let $n$ be a fixed natural number.

(i) By a known result
first observed by Euler (cf. [D99, p.\,260] and also [P]), there are $x,y,z\in\Z$ such that $2n+1=2x^2+4y^2+z^2$
and hence $n=x^2+2y^2+\lfloor z^2/2\rfloor$.

Suppose that $n\not=x^2+2y^2+\lfloor(3z)^2/3\rfloor =x^2+2y^2+3z^2$ for all $x,y,z\in\Z$.
Then $n$ is even by a known result (cf. [D39, p.\,112-113] or [P]). By [D39, p.\,112-113],
$$E(3x^2+6y^2+z^2)=\{3k+2:\ k\in\N\}\cup\{4^k(16l+14):\ k,l\in\N\}.$$
Since $3n+1$ is odd, for some $x,y,z\in\Z$ we have $3n+1=3x^2+6y^2+z^2$ and hence $n=x^2+2y^2+\lfloor z^2/3\rfloor$.

By [D39, p.\,112-113],
$$E(4x^2+8y^2+z^2)=\bigcup_{k\in\N}\{4k+2,4k+3\}\cup\{4^k(16l+14):\ k,l\in\N\}.$$
So there are $x,y,z\in\Z$ such that $4n+1=4x^2+8y^2+z^2$ and hence $n=x^2+2y^2+\lfloor z^2/4\rfloor$.

By [D39, p.\,112-113],
$$E(5x^2+10y^2+z^2)=\bigcup_{k,l\in\N}\{25^k(5l+2),25^k(5l+3)\}.$$
Thus, for some $x,y,z\in\Z$ we have $5n+1=5x^2+10y^2+z^2$ and hence $n=x^2+2y^2+\lfloor z^2/5\rfloor$.

(ii) By [D39, p.\,112-113],
$$E(3x^2+y^2+z^2)=\{9^k(9l+6):\ k,l\in\N\}.$$
So, there are $x,y,z\in\Z$ such that $3n+1=3x^2+(3y)^2+z^2$ and hence $n=x^2+3y^2+\lfloor z^2/3\rfloor$.

By [D39, p.\,112-113],
$$E(4x^2+12y^2+z^2)=\bigcup_{k\in\N}\{4k+2,4k+3\}\cup\{9^k(9l+6):\ k,l\in\N\}.$$
Choose $\da\in\{0,1\}$ such that $4n+\da\not\eq0\pmod 3$.
Then, for some $x,y,z\in\Z$ we have $4n+\da=4x^2+12y^2+z^2$ and hence $n=x^2+3y^2+\lfloor z^2/4\rfloor$.

If $6n+r=6x^2+18y^2+z^2$ for some $r\in\{0,1,3,4\}$ and $x,y,z\in\Z$, then $n=x^2+3y^2+\lfloor z^2/6\rfloor$.
Now suppose that $6n+r\not=6x^2+18y^2+z^2$ for any $r\in\{0,1,3,4\}$ and $x,y,z\in\Z$. By [D39, p.\,112-113],
$$S:=E(6x^2+18y^2+z^2)=\bigcup_{k\in\N}\{3k+2,9k+3\}\cup\{4^k(8l+5):\ k,l\in\N\}.$$
So $6n+1$ or $6n+4$ is congruent to $5$ modulo $8$. If $6n+4\eq5\pmod 8$, then $6n+1\eq2\pmod 8$ which
contradicts that $6n+1\in S$. So, $6n+1\eq 5\pmod 8$ and hence $6n+3\eq 7\pmod 8$. By $6n+3\in S$, we must have
$3\mid n$. As $6n\eq0\pmod 9$ and $6n\eq 4\pmod 8$, by $6n\in S$ we have $6n=4(8q+5)$ for some $q\in\Z$.
As $6n+4=4(8q+6)\not\in S$, we get a contradiction.

As conjectured by Sun [S07] and confirmed in [GPS], there are $x,y,z\in\Z$ such that $n=x^2+3y^2+T_z$ and hence $n=x^2+3y^2+\lfloor (2z+1)^2/8\rfloor$.

(iii) By [D39, p.\,112-113], $E(8x^2+40y^2+z^2)$ coincides with
$$\bigcup_{k\in\N}\{4k+2,4k+3,8k+5,32k+28\}\cup\bigcup_{k,l\in\N}\{25^k(25l+5),25^k(25l+20)\}.$$
Choose $\da\in\{0,1\}$ such that $8n+\da\not\eq0\pmod 5$. Then $8n+\da\not\in E(8x^2+40y^2+z^2)$.
So, for some $x,y,z\in\Z$ we have $8n+\da=8x^2+40y^2+z^2$ and hence $n=x^2+5y^2+\lfloor z^2/8\rfloor$.

By [D39, p.\,112-113],
$$E(4x^2+24y^2+z^2)=\bigcup_{k\in\N}\{4k+2,4k+3\}\cup\{9^k(9l+3):\ k,l\in\N\}.$$
Choose $\da\in\{0,1\}$ such that $4n+\da\not\eq0\pmod 3$. Then $4n+\da\not\in E(4x^2+24y^2+z^2)$.
Hence there are $x,y,z\in\Z$ such that $4n+\da=4x^2+24y^2+z^2$ and thus $n=x^2+6y^2+\lfloor z^2/4\rfloor$.

(iv) By [JP] or [D39, p.\,112-113], for some $x,y,z\in\Z$ we have $8n+1=16x^2+16y^2+z^2$ and hence $n=2x^2+2y^2+\lfloor z^2/8\rfloor$.

In view of [D39, p.\,112-113],
$$E(6x^2+9y^2+z^2)=\{3k+2:\, k\in\N\}\cup\{9^k(9l+3):\,k,l\in\N\}.$$
So, there are $x,y,z\in\Z$ such that $3n+1=6x^2+9y^2+z^2$ and hence $n=2x^2+3y^2+\lfloor z^2/3\rfloor$.

By [D39, p.\,112-113],
$$E(3x^2+3y^2+2z^2)=\{9^k(3l+1):\,k,l\in\N\}.$$
So there are $x,y,z\in\Z$ such that $6n+5=3x^2+3y^2+2z^2$.  Since $x\not\eq y\pmod2$, without loss of generality we may assume that
$2\mid x$ and $2\nmid y$. Thus
$$n=2\l(\f x2\r)^2+\f{y^2-1}2+\f{z^2-1}3=2\l(\f x2\r)^2+\l\lfloor\f{y^2}2\r\rfloor+\l\lfloor\f{z^2}3\r\rfloor.$$

\medskip

So far we have completed the proof of Theorem 1.2. \qed

\heading{3. Proofs of Theorems 1.3-1.4}\endheading

\medskip
\noindent{\it Proof of Theorem 1.3}. (i) Clearly, $0=\lceil0^2/m\rceil+\lceil0^2/m\rceil+\lceil0^2/m\rceil$,
$1=\lceil 1^2/3\rceil+\lceil0^2/3\rceil+\lceil0^2/3\rceil$ and $2=\lceil 1^2/3\rceil+\lceil1^2/3\rceil+\lceil0^2/3\rceil$.
for any $m\in\{2,3,4,5\}$.
So we just consider required representations for $n\in\{3,4,5,\ldots\}$.

If $n$ is even, then $2n-2\eq2\pmod 4$, hence by the Gauss-Legendre theorem there are integers $x,y,z$ with $2\nmid yz$
such that $2n-2=(2x)^2+y^2+z^2$ and thus
$$n=2x^2+\f{y^2+1}2+\f{z^2+1}2=\f{(2x)^2}2+\l\lceil\f{y^2}2\r\rceil+\l\lceil\f{z^2}2\r\rceil.$$
When $n\eq1\pmod 4$, we have $2n-1\eq1\pmod 8$
and hence by the Gauss-Legendre theorem there are $x,y,z\in\Z$ with $2\nmid z$ such that
$2n-1=(2x)^2+(2y)^2+z^2$ and thus
$$n=2x^2+2y^2+\f{z^2+1}2=\f{(2x)^2}2+\f{(2y)^2}2+\l\lceil\f{z^2}2\r\rceil.$$
If $n\eq3\pmod 4$, then $2n-3\eq3\pmod 8$, hence there are odd integers $x,y,z$ such that
$2n-3=x^2+y^2+z^2$ and thus
$$n=\f{x^2+1}2+\f{y^2+1}2+\f{z^2+1}2=\l\lceil\f{x^2}2\r\rceil+\l\lceil\f{y^2}2\r\rceil+\l\lceil\f{z^2}2\r\rceil.$$
This proves (1.7) for $m=2$.

Now we show (1.7) for $m=3$. Clearly we cannot have $\{3n-4,3n-6\}\se\{4^k(8l+7):\ k,l\in\N\}$ and hence either $3n-4$ or $3n-6$
can be written as the sum of three squares. If $3n-4=x^2+y^2+z^2$ for some $x,y,z\in\Z$, then exactly one of $x,y,z$
(say, $x$) is divisible by $3$, hence
$$n=3\l(\f x3\r)^2+\f{y^2+2}3+\f{z^2+2}3=\l\lceil\f{x^2}3\r\rceil+\l\lceil\f{y^2}3\r\rceil+\l\lceil\f{z^2}3\r\rceil.$$
When $3n-6=x^2+y^2+z^2$ with $x,y,z\in\Z$ not all zero, by [S16, Lemma 2.2] there are $u,v,w\in\Z$ with $3\nmid uvw$
such that $3n-6=u^2+v^2+w^2$ and hence
$$n=\f{u^2+2}3+\f{v^2+2}3+\f{w^2+2}3=\l\lceil\f{u^2}3\r\rceil+\l\lceil\f{v^2}3\r\rceil+\l\lceil\f{w^2}3\r\rceil.$$

As $4n-3\not\in\{4^k(8l+7):\ k,l\in\N\}$, by the Gauss-Legendre theorem there are $x,y,z\in\Z$ such that
$4n-3=(2x)^2+(2y)^2+(2z+1)^2$ and hence
$n=x^2+y^2+\lceil(2z+1)^2/4\rceil$. This proves (1.7) for $m=4$.

Now we prove (1.7) for $m=5$ by modifying our proof of the last equality in (1.5).
If $\{5n-5,5n-6,5n-9\}\se E:=\{4^k(8l+7):\ k,l\in\N\}$, then $5n-6\eq7\pmod 8$ and hence $5n-5\eq6\pmod8$
which leads a contradiction. So, for some $r\in\{0,1,4\}$ we can write $5n-5-r>5$ as the sum of three squares.
If $5n-5-r=m^2$ for some integer $m>2$ which is not a power of two, then by Lemma 2.1
we have $5(n-1)-r=x^2+y^2+z^2$ for some $x,y,z\in\Z$ with $x^2,y^2,z^2\not=5n-5-r$.
If $5(n-1)-r=(2^k)^2$ for some $k\in\Z^+$, then $r\in\{1,4\}$, and $5(n-1)-(5-r)=4^k+2r-5\eq2r-5\eq \pm3\pmod{8}$
and hence $5(n-1)-(5-r)\not\in E$. So, for a suitable choice of $r\in\{0,1,4\}$,
we can write $5(n-1)-r=x^2+y^2+z^2$ with $x,y,z\in\Z$ and $x^2,y^2,z^2\not=5(n-1)-r$. Clearly, one of $x^2,y^2,z^2$, say $z^2$, is congruent to $-r$ modulo $5$.
Then $x^2+y^2$ is a positive multiple of $5$. By Lemma 2.2, $x^2+y^2=\bar x^2+\bar y^2$ for some $\bar x,\bar y\in\Z$ with $5\nmid \bar x\bar y$.
Without loss of generality, we may assume that $\bar x^2\eq 1\pmod 5$ and $\bar y^2\eq4\pmod 5$. Therefore,
$$n=\f{\bar x^2+4}5+\f{\bar y^2+1}5+\f{z^2+r}5=\l\lceil\f{\bar x^2}5\r\rceil+\l\lceil\f{\bar y^2}5\r\rceil+\l\lceil\f{z^2}5\r\rceil.$$

Now we show (1.7) for $m=6$. If $n$ is odd, then $6n-9\eq1\pmod4$, hence by the Gauss-Legendre theorem and [S16, Lemma 2.2] we can write
$6n-9=x^2+y^2+z^2$ with $x,y,z\in\Z$, $2\nmid x$, $2\mid y$, $2\mid z$  and $3\nmid xyz$, therefore
$$n=\f{x^2+5}6+\f{y^2+2}6+\f{z^2+2}6=\l\lceil\f{x^2}6\r\rceil+\l\lceil\f{y^2}6\r\rceil+\l\lceil\f{z^2}6\r\rceil.$$
Now assume that $n$ is even. Then $6n-10\eq2\pmod{12}$. By the Gauss-Legendre theorem we can write $6n-10=x^2+y^2+z^2$
with $x,y,z\in\Z$, $2\nmid xy$ and $2\mid z$. Note that exactly one of $x,y,z$ is divisible by $3$.
If $3\nmid xy$ and $3\mid z$, then $x^2\eq y^2\eq1\pmod 6$ and $z^2\eq0\pmod 6$, hence
$$n=\f{x^2+5}6+\f{y^2+5}6+\f{z^2}6=\l\lceil\f{x^2}6\r\rceil+\l\lceil\f{y^2}6\r\rceil+\l\lceil\f{z^2}6\r\rceil.$$
If $3\nmid z$, then exactly one of $x$ and $y$, say $x$, is divisible by 3, hence
$x^2\eq 3\pmod 6$, $y^2\eq1\pmod 6$ and $z^2\eq 4\pmod 6$, and thus
$$n=\f{x^2+3}6+\f{y^2+5}6+\f{z^2+2}6=\l\lceil\f{x^2}6\r\rceil+\l\lceil\f{y^2}6\r\rceil+\l\lceil\f{z^2}6\r\rceil.$$

Now we prove (1.7) for $m=15$. By the proof of the last equality in (1.5) for $m=15$, for a suitable choice of $r\in\{1,4,10\}$ we have
$15(n-3)+5+r=x^2+y^2+z^2$ for some $x,y,z\in\Z$ with $x^2\eq1\pmod{15}$, $y^2\eq4\pmod{15}$ and $z^2\eq r\pmod{15}$.
It follows that
$$n=\f{x^2+14}{15}+\f{y^2+11}{15}+\f{z^2+15-r}{15}=\l\lceil\f{x^2}{15}\r\rceil+\l\lceil\f{y^2}{15}\r\rceil+\l\lceil\f{z^2}{15}\r\rceil.$$

(ii) Now we turn to prove (1.8). Apparently, $0=0^2+3\times0^2+\lceil 0^2/2\rceil$. Let $n\in\Z^+$.
If $2n-1\eq 5\pmod 8$ then $4\nmid 2n$. So, we may choose $\da\in\{0,1\}$ such that
$2n-\da\not\in\{4^k(8l+5):\ k,l\in\N\}$. By [D39, p.\,112-113],
$$E(2x^2+6y^2+z^2)=\{4^k(8l+5):\ k,l\in\N\}.$$
So there are $x,y,z\in\Z$ such that $2n-\da=2x^2+6y^2+z^2$ and hence
$n=x^2+3y^2+\lceil z^2/2\rceil$.

Obviously, $0=0^2+3\times0^2+\lceil0^2/10\rceil$. Let $n\in\Z^+$.
By [D39, p.\,112-113],
$$\align T:=E(10x^2+30y^2+z^2)=\bigcup_{k,l\in\N}\{4^k(8l+5),9^k(9l+6),25^k(5l+2),25^k(5l+3)\}.
\endalign$$
If $10n-r\not\in T$ for some $r\in\{0,1,4,5,6,9\}$, then there are $x,y,z\in\Z$ such that $10n-r=10x^2+30y^2+z^2$
and hence $n=x^2+3y^2+\lceil z^2/10\rceil$. Now we suppose that $10n-r\in T$ for all $r=0,1,4,5,6,9$ and want to deduce a contradiction.
If $3\mid n(n+1)$, then by $10n-1\in T$ we have $10n-1\eq 5\pmod8$ and hence $10n-4\eq2\pmod8$ which contradicts $10n-4\in T$.
When $n\eq1\pmod3$, by $10n-9\in T$ we must have $10n-9\eq5\pmod 8$ and thus $10n\eq6\pmod 8$,
hence $10n\eq0\not\eq5\pmod{25}$  by $10n\in T$, and thus by $10n-5\in T$ we have
$10n-5\eq5\pmod 8$ which contradicts $10n\eq6\pmod 8$.

(iii) Choose $\da\in\{0,1\}$ with $n\eq\da\pmod2$. Then $12n+5-4\da\not\eq0\pmod 3$. By [D39, pp.\,112-113],
$$E(3x^2+y^2+z^2)=\{9^k(9l+6):\ k,l\in\N\}.$$
So, there are $u,v,w\in\Z$ such that $12n+5-4\da=3u^2+v^2+w^2$.
If $v$ and $w$ are both even, then $5\eq3u^2\pmod 4$ which is impossible.
Without loss of generality, we assume that $w=2z+1$ with $z\in\Z$. Then
$$3u^2+v^2\eq 12n+5-4\da-1\eq 4\pmod 8.$$
Hence, by [S15a, Lemma 3.2] we can write $3u^2+v^2$ as $3(2x+1)^2+(2y+1)^2$ with $x,y\in\Z$.
Therefore,
$$12n+5-4\da=3(2x+1)^2+(2y+1)^2+(2z+1)^2=12x(x+1)+4y(y+1)+4z(z+1)+5$$
and hence
$$3n-\da=3x(x+1)+y(y+1)+z(z+1).$$
Note that $m(m+1)\not\eq1\pmod3$ for any $m\in\Z$.
If $y(y+1),z(z+1)\not\eq0\pmod 3$, then $-\da\eq 2+2\pmod3$ which is impossible.
Without loss of generality we assume that $3\mid y(y+1)$. Then
$$n=x(x+1)+\f{y(y+1)}3+\f{z(z+1)+\da}3=x(x+1)+\f{y(y+1)}3+\l\lceil\f{z(z+1)}3\r\rceil.$$

Let $\da\in\{0,1\}$ with $n\eq\da\pmod2$. Then $12n+3-4\da$ is congruent to $0$ or $2$ modulo $3$.
As $12n+3-4\da\eq 3\pmod 8$, there are odd integers $u,v,w$ such that $12n+3-4\da=u^2+v^2+w^2$.
If $\da=0$, then by [S16, Lemma 2.2] we can write $u^2+v^2+w^2$ as $r^2+s^2+t^2$ with $r,s,t\in\Z$ and $\gcd(rst,6)=1$.
So, there are $x,y,z\in\Z$ such that
$$12n+3-4\da=(2x+1)^2+(2y+1)^2+(2z+1)^2=4x(x+1)+4y(y+1)+4z(z+1)+3$$
and $2x+1,2y+1\not\eq0\pmod3$.  As $x,y\not\eq1\pmod3$, both $x(x+1)$ and $y(y+1)$ are divisible by $3$. Thus
$$n=\f{x(x+1)}3+\f{y(y+1)}3+\f{z(z+1)+\da}3=\l\lceil\f{x(x+1)}3\r\rceil+\l\lceil\f{y(y+1)}3\r\rceil+\l\lceil\f{z(z+1)}3\r\rceil.$$
Note that $\{m(m+1)/3:\ m\in\Z\ \&\ 3\mid m(m+1)\}=\{q(3q+1):\ q\in\Z\}$.

The proof of Theorem 1.3 is now complete. \qed

\medskip
\noindent{\it Proof of Theorem 1.4}. (i) Let $n\in\N$. If $2n+1\in\{4^k(8l+7):\ k,l\in\N\}$, then $2n\eq6\pmod 8$ and hence $2n\not\in\{4^k(8l+7):\ k,l\in\N\}$.
So, for some $\da\in\{0,1\}$ we have $2n+\da\not\{4^k(8l+7):\ k,l\in\N\}$, and hence by the Gauss-Legendre theorem there are $x,y,z\in\Z$
such that $2n+\da=x^2+y^2+z^2$ and hence $n=\lfloor (x^2+y^2+z^2)/2\rfloor$. Note also that $n=T_x+T_y+T_z$ for some $x,y,z\in\Z$.
This proves the desired result for $a=2$.

Now we handle the case $a>2$. Clearly, for some $r\in\{0,2\}$ we have $an+r\not\in\{4^k(8l+7):\ k,l\in\N\}$,
hence for some $x,y,z\in\Z$ we have $an+r=x^2+y^2+z^2$ and thus $n=\lfloor(x^2+y^2+z^2)/a\rfloor$.
Take $\da\in\{0,1\}$ with $an\eq \da\pmod2$. Then, there exist $x,y,z\in\Z$ such that $(an+da)/2=T_x+T_y+T_z$
and hence $n=\lfloor(x(x+1)+y(y+1)+z(z+1))/a\rfloor$.

(ii) Suppose that $a$ is odd. As $16n+3a^2\eq 3\pmod 8$, by the Gauss-Legendre symbol $16n+3a^2$ can be expressed as the sum of three odd squares.
For any odd integer $w$, either $w$ or $-w$ is congruent to $a$ modulo $4$. Thus, there are $x,y,z\in\Z$ such that
$$16n+3a^2=(4x+a)^2+(4y+a)^2+(4z+a)^2,\ \t{i.e.},\ 2n=2(x^2+y^2+z^2)+a(x+y+z).$$
Hence $n=x^2+y^2+z^2+\lfloor \f a2(x+y+z)\rfloor$ as desired.

Now assume that $\gcd(a,6)=1$. Choose $\da\in\{0,1\}$ such that $n\eq\da\pmod 2$. As $12(3n+\da)+3a^2\eq3\pmod 8$, there are odd integers $u,v,w$
such that $12(3n+\da)+3a^2=u^2+v^2+w^2$.  Applying [S16, Lemma 2.2], we can write $u^2+v^2+w^2$
as $r^2+s^2+t^2$, where $r,s,t$ are integers with
$$r\eq u_0\eq u\eq1\pmod 2,\ s\eq v\eq 1\pmod 2,\ t\eq w\eq1\pmod2,\ \t{and}\ 3\nmid rst.$$
Thus $r$ or $-r$ has the form $6x+a$, $s$ or $-s$ has the form $6y+a$, and $t$ or $-t$ has the form $6z+a$, where $x,y,z\in\Z$.
Therefore,
$$\align 12(3n+\da)+3a^2=&(6x+a)^2+(6y+a)^2+(6z+a)^2
\\=&12(3x^2+ax+3y^2+ay+3z^2+3z)+3a^2
\endalign$$
and hence
$$n=x^2+y^2+z^2+\f {a(x+y+z)-\da}3=x^2+y^2+z^2+\l\lfloor\f a3(x+y+z)\r\rfloor.$$

Now we suppose that $2\mid a$ and $3\nmid a$. If $9n+3(a/2)^2+3r\in \{4^k(8l+7):\ k,l\in\N\}$ for all $r=1,2,3$, then $9n+3(a/2)^2+6\eq7\pmod 8$ and hence $9n+3(a/2)^2+9\eq 2\pmod 8$
which leads a contradiction. So, by the Gauss-Legendre theorem, for some $r\in\{1,2,3\}$ and $u,v,w\in\Z$ we have
$9n+3(a/2)^2+3r=u^2+v^2+w^2$. By [S16, Lemma 2.2] we can write $9n+3(a/2)^2+3r=\bar u^2+\bar v^2+\bar w^2$, where $\bar u,\bar v,\bar w\in\Z$ and $3\nmid\bar u\bar v\bar w$.
So there are $x,y,z\in\Z$ such that
$$9n+3r+3\l(\f a2\r)^2=\l(3x+\f a2\r)^2+\l(3y+\f a2\r)^2+\l(3z+\f a2\r)^2,$$
i.e.,
$$3n+r-1=x(3x+a)+y(3y+a)+z(3z+a).$$
It follows that
$$n=x^2+y^2+z^2+\f{a(x+y+z)-(r-1)}3=x^2+y^2+z^2+\l\lfloor\f a3(x+y+z)\r\rfloor.$$

(iii) Obviously, $0=p_8(0)/2+\lceil p_8(0)/2\rceil+\lceil p_8(0)/2\rceil$. Now we let $n>0$
and choose $\da\in\{0,1\}$ with $n\not\eq\da\pmod2$. As $6n-3\da$ is congruent to $1$ or $2$ modulo $4$,
by the Gauss-Legendre theorem we can write $6n-3\da$ as the sum of three squares
and hence by [S16, Lemma 2.2] there are $x,y,z\in\Z$ such that
$$6n-3\da=(3x-1)^2+(3y-1)^2+(3z-1)=3p_8(x)+1+(3p_8(y)+1)+(3p_8(z)+1).$$
Clearly, $3x-1,3y-1,3z-1$ cannot be all odd or all even.
Without loss of generality, we may assume that
$$3x-1\eq1\pmod2,\ 3y-1\eq0\pmod2\ \t{and}\ 3z-1\eq1-\da\eq n\pmod2.$$
Then $p_8(x)=((3x-1)^2-1)/3$ is even, $p_8(y)$ is odd, and $p_8(z)\eq -\da\pmod2$. Therefore
$$n=\f{p_8(x)}2+\f{p_8(y)+1}2+\f{p_8(z)+\da}2=\f{p_8(x)}2+\l\lceil\f{p_8(y)}2\r\rceil+\l\lceil\f{p_8(z)}2\r\rceil.$$
This concludes our proof. \qed

\heading{4. Proof of Theorem 1.5}\endheading

For $a,b,c,n\in\Z^+$, define
$$r_{(a,b,c)}(n)=|\{(x,y,z)\in\Z^3:\ ax^2+by^2+cz^2=n\}|\tag4.1$$
and
$$H_{(a,b,c)}(n):=\prod_{p\nmid 2abc}\l(\f{p^{\ord_p(n)+1}-1}{p-1}-\l(\f{-abc}p\r)\f{p^{\ord_p(n)}-1}{p-1}\r),\tag4.2$$
where $(\f{\cdot}p)$ is the Legendre symbol. Clearly,
$$H_{(a,b,c)}(n)\gs\prod_{p\nmid 2abc}\f{p^{\ord_p(n)+1}-1-(p^{\ord_p(n)}-1)}{p-1}=\prod_{p\nmid2abc}p^{\ord_p(n)}.\tag4.3$$
In 1907 Hurwitz (cf. [D99, p.\,271]) showed that $r_{(1,1,1)}(n^2)=6H_{(1,1,1)}(n)$ which is just (2.1).
In 2013 S. Cooper and H. Y. Lam [CL] deduced some similar formulas for
$$r_{(1,1,2)}(n^2),\ r_{(1,1,3)}(n^2),\ r_{(1,2,2)}(n^2),\ r_{(1,3,3)}(n^2).$$

\proclaim{Lemma 4.1} For any integer $n>1$, there are $x,y,z\in\Z$ with $|x|<n$ and $|y|<n$ such that $x^2+y^2+2z^2=n^2$.
\endproclaim
\Proof. By Cooper and Lam [CL, Theorem 1.2],
$$\aligned r_{(1,1,2)}(n^2)=\cases 4H_{(1,1,2)}(n)&\t{if}\ 2\nmid n,
\\12H_{(1,1,2)}(n)&\t{if}\ 2\mid n.\endcases
\endaligned\tag4.4$$
If $n$ is odd, then there is an odd prime $p$ dividing $n$, hence $r_{(1,1,2)}(n^2)=4H_{(1,1,2)}(n)>4$ with the help of (4.3).
If $n$ is even, then $r_{(1,1,2)}(n^2)=12H_{(1,1,2)}(n)\gs 12>4$. Clearly,  $x^2+y^2+2z^2=n^2$
for $(x,y,z)=(\pm n,0,0),(0,\pm n,0)$. So, there are $x,y,z\in\Z$ with $x^2,y^2\not=n^2$ such that $x^2+y^2+2z^2=n^2$.
This concludes the proof. \qed

\proclaim{Lemma 4.2} Suppose that $n\in\Z^+$ is not a power of two. Then there are $x,y,z\in\Z$ with $|x|<n$ and $|y|<n$ such that $x^2+y^2+5z^2=n^2$.
\endproclaim
\Proof. As conjectured by Cooper and Lam [CL, Conjecture 8.1] and proved by Guo et al. [GPQ],
$$r_{(1,1,5)}(n^2)=2(5^{\ord_5(n)+1}-3)H_{(1,1,5)}(n).\tag4.5$$
If $5\mid n$, then $2(5^{\ord_5(n)+1}-3)>4$. If $n$ has a prime divisor $p\not=2,5$, then
$H_{(1,1,5)}(n)>1$ by (4.3). Since $n>1$ is not a power of two, we have $r_{(1,1,5)}(n^2)>4$.
Clearly,  $x^2+y^2+5z^2=n^2$
for $(x,y,z)=(\pm n,0,0),(0,\pm n,0)$. So, there are $x,y,z\in\Z$ with $x^2,y^2\not=n^2$ such that $x^2+y^2+5z^2=n^2$.
This ends the proof. \qed

\Remark\ 4.1. Note that Lemmas 4.1 and 4.2 are similar to Lemma 2.1.

\medskip
\noindent{\it Proof of Theorem 1.5}. (i) Let $n\in\N$. Clearly, $n=p_8(x)+p_8(y)+2p_8(z)$ if and only if $3n+4=(3x-1)^2+(3y-1)^2+2(3z-1)^2$.
In view of [D39, pp.\,112-113],
$$E(x^2+y^2+2z^2)=\{4^k(16l+14):\ k,l\in\N\}.$$
If $3n+4=4^k(16l+14)$ for some $k,l\in\N$, then for some $q\in\Z^+$
we have $l=3q-1$ and hence
$$n=\f{4^k(16(3q-1)+14)-4}3=4^{k+2}q-\f23(4^k+2).$$
If $n$ has the form $4^{k+2}q-\f23(4^k+2)$ with $k\in\N$ and $q\in\Z^+$,
then $n\not=p_8(x)+p_8(y)+2p_8(z)$ for all $x,y,z\in\Z$.

Now assume that $n$ is not of the form $4^{k+2}q-\f23(4^k+2)$ with $k\in\N$ and $q\in\Z^+$.
Then there are $r,s,t\in\Z$ such that $3n+4=r^2+s^2+2t^2$.
In view of Lemma 4.1, we may assume that $r^2,s^2\not=3n+4$.
Clearly $r$ and $s$ cannot be both divisible by $3$. Without loss of generality, we assume that
$3\nmid r$. As $s^2+2t^2=3n+4-r^2$ is a positive multiple of $3$, by [S15a, Lemma 2.1] we can rewrite it as
$u^2+2v^2$ with $u,v\in\Z$ and $3\nmid uv$. Thus there are $x,y,z\in\Z$ such that
$$\align 3n+4=&r^2+u^2+2v^2=(3x-1)^2+(3y-1)^2+2(3z-1)^2
\\=&3p_8(x)+1+(3p_8(y)+1)+2(3p_8(z)+1)
\endalign$$
and hence $n=p_8(x)+p_8(y)+2p_8(z)$.

By the above, there are $x,y,z\in\Z$ with $2n+1=p_8(x)+p_8(y)+2p_8(z)$.
Without loss of generality, we may assume that $p_8(x)$ is even and $p_8(y)=y(3y-2)$ is odd.
Clearly, $w=(1-y)/2\in\Z$ and $p_8(y)-1=4p_8(w)$. So, $2n=p_8(x)+2p_8(z)+4p_8(w)$.
Note also that
$$n=\f{p_8(x)}2+\f{p_8(y)-1}2+p_8(z)=\l\lfloor\f{p_8(x)}2\r\rfloor+\l\lfloor\f{p_8(y)}2\r\rfloor+p_8(z).$$
Therefore (1.13) and (1.14) hold.

(ii) Fix a nonnegative integer $n$. If $6n+5\eq7\pmod 8$, then $6n+8\eq2\pmod 8$.
So, for a suitable choice of $\da\in\{0,1\}$ we have $6n+5+3\da\not\in E(x^2+y^2+z^2)=\{4^k(8l+7):\ k,l\in\N\}$
and hence $6n+5+3\da=u^2+v^2+w^2$ for some $u,v,w\in\Z$. Two of $u,v,w$ have the same parity.
Without loss of generality, we assume that $u+v=2s$ and $u-v=2t$ for some $s,t\in\Z$.
Hence $6n+5+3\da=w^2+2s^2+2t^2$. If $(6n+5+3\da)=2m^2$ for some $m\in\Z^+$, then by Lemma 4.1 there
are $r,s_1,t_1\in\Z$ with $s_1^2,t_1^2\not=m^2$ such that $m^2=s_1^2+t_1^2+2r^2$
and hence $6n+5+3\da=(2r)^2+2s_1^2+2t_1^2$ with $2s_1^2,2t_1^2\not=6n+5+3\da$.
So, we may simply suppose that $6n+5+3\da=w^2+2s^2+2t^2$ with $2s^2,2t^2\not=6n+5+3\da$.
Clearly, one of $s$ and $t$ is not divisible by $3$. Without loss of generality we assume that $t^2=(3x-1)^2$ with $x\in\Z$ As $w^2+2s^2=6n+5+3\da-2t^2$
is a positive multiple of $3$, by [S15a, Lemma 2.1] we can write $w^2+2s^2$ as $(3z-1)^2+2(3y-1)^2$ with $y,z\in\Z$. Thus
$$6n+5+3\da=(3z-1)^2+2(3y-1)^2+2(3x-1)^2=3p_8(z)+1+2(3p_8(x)+3p_8(y)+2)$$
and hence
$$n=p_8(x)+p_8(y)+\f{p_8(z)-\da}2=p_8(x)+p_8(y)+\l\lfloor\f{p_8(z)}2\r\rfloor.$$
This proves (1.15).  In view of (1.12), both (1.16) and (1.17) follow from (1.15).

(iii) Let $n\in\N$. As $12n+9\eq1\pmod 4$, by the Gauss-Legendre theorem we can write $12n+9$
as the sum of three squares. In view of [S16, Lemma 2.2], there are $u,v,w\in\Z$ with $3\nmid uvw$
such that $12n+9=u^2+v^2+w^2$. Clearly, exactly one of $u,v,w$ is odd. Without loss of generality
we may assume that $u=2(3x-1)$, $v=2(3y-1)$ and $w=3z-1$ with $x,y,z\in\Z$. Thus
$$12n+9=4(3x-1)^2+4(3y-1)^2+(3z-1)^2=12p_8(x)+12p_8(y)+3p_8(z)+9$$
and hence (1.18) follows.

(iv) By Dickson [D39, pp.\,112-113],
$$E(5x^2+5y^2+z^2)=\bigcup_{k\in\N}\{5k+2,\,5k+3\}\cup\{4^k(8l+7):\,k,l\in\N\}.$$
If $15n+11+3r$ belongs to this set for all $r=0,1,3$, then
$15n+11$ is odd, hence $15n+11\eq7\pmod 8$ and $15n+11+3\eq2\pmod 8$
which leads a contradiction. So, there is a choice of $r\in\{0,1,3\}$ such that
$15n+11+3r\not\in E(5x^2+5y^2+z^2)$. Hence, for some $u,v,w\in\Z$ we have $15n+11+3r=5u^2+5v^2+w^2$.
If $15n+11+3r=5m^2$ for some positive integer $m$ which is not a power of two, then
by Lemma 4.2 there are $u_1,v_1,w_1\in\Z$ with $u_1^2,v_1^2\not=m^2$ such that $m^2=u_1^2+v_1^2+5w_1^2$
and hence $15n+11+3r=5u_1^2+5v_1^2+(5w_1)^2$ with $5u_1^2,5v_1^2\not=15n+11+3r$.
If $15n+11+3r=5\times 2^a$ for some $a\in\N$, then $a\gs2$, $r=3$, $15n+11+3\times1=5\times2^a-6\eq2\pmod 4$
and hence $15n+11+3\not\in E(5x^2+5y^2+z^2)$. So, we may simply assume that $15n+11+3r=5u^2+5v^2+w^2$
with $5u^2,5v^2<15n+11+3r$. Clearly, $u$ or $v$ is not divisible by $3$. Without loss of generality we suppose that
$u^2=(3x-1)^2$ for some $x\in\Z$. As $5v^2+w^2=15n+11+3r-5u^2>0$ is a positive multiple of $3$, by [S15a, Lemma 2.1]
we can write $5v^2+w^2$ as $5(3y-1)^2+(3z-1)^2$ with $y,z\in\Z$. Thus
$$\align 15n+11+3r=&5(3x-1)^2+5(3y-1)^2+(3z-1)^2
\\=&5(3p_8(x)+1)+5(3p_8(y)+1)+3p_8(z)+1
\endalign$$
and hence
$$n=p_8(x)+p_8(y)+\f{p_8(z)-r}5=p_8(x)+p_8(y)+\l\lfloor\f{p_8(z)}5\r\rfloor.$$
This proves (1.19). In view of (1.12), (1.20) follows from (1.19).

The proof of Theorem 1.5 is now complete. \qed

\heading{5. Some further conjectures}\endheading

\proclaim{Conjecture 5.1} For any $n\in\N$, there are $x,y,z\in\N$ such that $8n+3=x^2+y^2+z^2$ and $x\eq1,3\pmod8$.
Also, for any $n\in\N$ with $n\not=20$, there are $x,y,z\in\Z$ with $x\eq\pm3\pmod 8$ such that $x^2+y^2+z^2=8n+3$.
\endproclaim
\Remark\ 5.1. In [S15a] the author conjectured that any $n\in\N$ can be written as the sum of two triangular numbers and a hexagonal number,
equivalently, $8n+3=(4x-1)^2+y^2+z^2$ for some $x,y,z\in\N$.

\proclaim{Conjecture 5.2} Let $a>2$ be an integer with $a\not=4,6$. Then any positive integer can be written as the sum of three elements of
the set $\{\lfloor x^2/a\rfloor:\,x\in\Z\}$ one of which is odd.
\endproclaim
\Remark\ 5.2. This is a refinement of Farhi's conjecture for $a\not=4,6$.

\proclaim{Conjecture 5.3} Let
$$T:=\l\{x^2+\l\lfloor \f x2\r\rfloor:\ x\in\Z\r\}=\l\{\l\lfloor\f{k(k+1)}4\r\rfloor:\ k\in\N\r\}.$$
Then each $n=2,3,4,\ldots$ can be expressed as $r+s+t$, where $r,s,t$ are elements of $T$
with $r\ls s\ls t$ and $2\nmid s$. Also, for any ordered pair $(b,c)$ among
$$(1,2),\,(1,3),\,(1,4),\,(1,5),\,(1,6),\,(1,8),\,(1,9),\,(2,2),\,(2,3),$$
each $n\in\N$ can be written as $x+by+cz$ with $x,y,z\in T$.
\endproclaim
\Remark\ 5.3. It is easy to see that $\{T_x:\ x\in\Z\}=\{p_6(-x)=x(2x+1):\ x\in\Z\}$.

\proclaim{Conjecture 5.4} {\rm (i) } Let $\al$ be a positive real number with $\al\not=1$ and $\al\ls 1.5$. Define
$$S(\al):=\{x^2+\lfloor \al x\rfloor:\ x\in\Z\}.$$
Then any positive integer can be written as the sum of three elements of $S(\al)$ one of which is odd.

{\rm (ii)} Let $0<\al\ls \beta\ls\gamma\ls 1.5$ such that two of $\al,\beta,\gamma$ are different from $1$ or $\{\al,\beta,\gamma\}=\{1,1/m\}$ for some $m=2,3,4,\ldots$. Then any $n\in\N$ can be written as
 $x^2+y^2+z^2+\lfloor \al x\rfloor+\lfloor \beta y\rfloor+\lfloor \gamma z\rfloor$ with $x,y,z\in\Z$. In particular, if $a,b,c\in\Z^+$ are not all equal to one, then
 $$\l\{x^2+y^2+z^2+\l\lfloor\f xa\r\rfloor+\l\lfloor\f yb\r\rfloor+\l\lfloor\f zc\r\rfloor:\ x,y,z\in\Z\r\}=\N.$$
\endproclaim
\Remark\ 5.4. Note that $2$ cannot be written as the sum of three elements of $S(11/4)$, and $4$ cannot be written as the sum of
three elements of $S(8/5)$ one of which is odd.
\medskip

\proclaim{Conjecture 5.5} Any integer $n>1$ can be written as $p+\lfloor k(k+1)/4\rfloor$, where $p$ is a prime and $k$ is a positive integer.
\endproclaim
\Remark\ 5.5. The author [S09] conjectured that $216$ is the only natural number not representable by $p+T_x$, where $p$ is prime or zero, and $x$ is an integer.
\medskip

Motivated by Theorem 1.5, we pose the following conjecture.

\proclaim{Conjecture 5.6} Let $a,b,c\in\Z^+$. Then
$$\l\{\l\lfloor\f{p_5(x)}a\r\rfloor+\l\lfloor\f{p_5(y)}b\r\rfloor+\l\lfloor\f{p_5(z)}c\r\rfloor:\ x,y,z\in\Z\r\}=\N.$$
When $(a,b,c)\not=(1,1,1),(1,1,2),(2,2,2)$, we have
$$\l\{\l\lfloor\f{p_7(x)}a\r\rfloor+\l\lfloor\f{p_7(y)}b\r\rfloor+\l\lfloor\f{p_7(z)}c\r\rfloor:\ x,y,z\in\Z\r\}=\N.$$
If $(a,b,c)\not=(1,1,1),(2,2,2)$, then
$$\l\{\l\lfloor\f{p_8(x)}a\r\rfloor+\l\lfloor\f{p_8(y)}b\r\rfloor+\l\lfloor\f{p_8(z)}c\r\rfloor:\ x,y,z\in\Z\r\}=\N.$$
\endproclaim

Now we present a general conjecture related to Theorems 1.1-1.3.

\proclaim{Conjecture 5.7} {\rm (i)} Let $a$ and $b$ be positive integers. If $c\in\Z^+$ is large enough, then
$$\l\{ax^2+by^2+\l\lfloor\f{z^2}c\r\rfloor:\ x,y,z\in\Z\r\}=\l\{ax^2+by^2+\l\lceil\f{z^2}c\r\rceil:\ x,y,z\in\Z\r\}=\N.$$
Also, for any sufficiently large $c\in\Z^+$ we have
$$\l\{ax^2+by^2+\l\lfloor\f{z(z+1)}c\r\rfloor:\ x,y,z\in\Z\r\}=\N$$
and
$$\l\{ax^2+by^2+\l\lceil\f{z(z+1)}c\r\rceil:\ x,y,z\in\Z\r\}=\N.$$

{\rm (ii)} For $a,b,c\in\Z^+$ with $2a\ls b+c$, if $(a,b,c)\not=(1,1,1),(3,3,3),(4,2,6)$ then
$$\l\{ax^2+\l\lfloor\f{y^2}b\r\rfloor+\l\lfloor\f{z^2}c\r\rfloor:\ x,y,z\in\Z\r\}=\N.$$
\endproclaim

\medskip

For $a,b\in\Z^+$, we define
$$\align S^*(a,b):=&\l\{c\in\Z^+:\ \l\{ax^2+by^2+\l\lceil\f{z^2}c\r\rceil:\, x,y,z\in\Z\r\}\not=\N\r\},
\\S_*(a,b):=&\l\{c\in\Z^+:\ \l\{ax^2+by^2+\l\lfloor\f{z^2}c\r\rfloor:\, x,y,z\in\Z\r\}\not=\N\r\},
\\T^*(a,b):=&\l\{c\in\Z^+:\ \l\{ax^2+by^2+\l\lceil\f{z(z+1)}c\r\rceil:\,x,y,z\in\Z\r\}\not=\N\r\},
\\T_*(a,b):=&\l\{c\in\Z^+:\ \l\{ax^2+by^2+\l\lfloor\f{z(z+1)}c\r\rfloor:\,x,y,z\in\Z\r\}\not=\N\r\}.
\endalign$$
Based on our computation we conjecture that
$$\gather S^*(1,1)=\{1,2,5\},\ S^*(1,2)=\{1,3\},\ S^*(1,3)=\{1,4\},\ S^*(1,4)=\{1,2,3,5\},
\\S^*(1,5)=\{1,2,3,5\},\ S^*(1,6)=\{1,2,3,4\},\ S^*(1,7)=\{1,2,4,8\},
\\S^*(1,8)=\{1,\ldots,6,9\},\ S^*(1,9)=\{1,\ldots,6\},\ S^*(1,10)=\{1,\ldots,6,8,12\},
\\S^*(2,2)=\{1,\ldots,5,9,10\},\ S^*(2,3)=\{1,2,8\};
\\S_*(1,2)=\{1\},\ S_*(1,3)=\{1,2,10\},\ S_*(1,4)=\{1,2,3,5\},\ S_*(1,5)=\{1,2,3,4,5\},
\\S_*(1,6)=\{1,3\},\ S_*(1,7)=\{1,2,3,4,5\},\ S_*(1,8)=\{1,2,3,5,9\},
\\S_*(1,9)=\{1,2,3,4,5,7\},\ S_*(1,10)=\{1,2,3,4,12\},
\\ S_*(1,11)=\{1,2,3,4,5,6,9\},\ S_*(1,12)=\{1,2,3,4,5,6,10\}
\\S_*(2,2)=\{1,2,3,4,5,6,10\},\ S_*(2,3)=\{1,2,8\},
\\ S_*(2,4)=\{1,2,5,6\},\ S_*(2,5)=\{1,2,3,5\};
\endgather$$
$$\gather T^*(1,1)=T^*(1,2)=\em,\ T^*(1,3)=1,\ T^*(1,4)=\{3\},\ T^*(1,5)=T^*(1,6)=\{1,2\},
\\T^*(1,7)=\{1,2,4\},\ T^*(1,8)=\{1\},\ T^*(1,9)=T^*(1,10)=T^*(1,11)=\{1,2,3\},
\\ T^*(2,2)=\{1,3\},\ T^*(2,3)=\{1,2\},\ T^*(2,4)=\{1,2,3\},\ T^*(3,4)=\{1,2,3\};
\\T_*(1,2)=\em,\ T_*(1,3)=\{1\},\ T_*(1,5)=\{1,2,3\},\ T_*(1,6)=\{1,2\},
\\T_*(1,7)=\{1,2,4\},\ T_*(1,8)=\{1\},\ T_*(1,10)=T_*(2,3)=\{1,2,3\}.
\endgather$$
Also, our computation suggests that
$$\l\{4x^2+4y^2+\l\lfloor\f{z^2}c\r\rfloor:\ x,y,z\in\Z\r\}=\N$$
for any integer $c>42$, and that
$$\l\{4x^2+4y^2+\l\lfloor\f{z(z+1)}c\r\rfloor:\ x,y,z\in\Z\r\}=\N$$
for any integer $c>27$. Note that $179\not=4x^2+4y^2+\lfloor z^2/42\rfloor$ for any $x,y,z\in\Z$
and that $29\not=4x^2+4y^2+\lfloor z(z+1)/27\rfloor$ for all $x,y,z\in\Z$.
\medskip

Motivated by Theorem 1.4(i), we pose the following conjecture similar to Conjecture 1.1.
\proclaim{Conjecture 5.8} Let $a,b,c$ be positive integers with $a\ls b\ls c$. If $c>1$, then
$$\l\{\l\lfloor\f{x^2}a+\f{y^2}b+\f{z^2}c\r\rfloor:\ x,y,z\in\Z\r\}=\N$$
If $(a,b,c)\not=(1,1,1),(1,1,3),(1,1,7),(1,3,3)$, then
$$\l\{\l\lfloor\f{x(x+1)}a+\f{y(y+1)}b+\f{z(z+1)}c\r\rfloor:\ x,y,z\in\Z\r\}=\N.$$
\endproclaim

\proclaim{Conjecture 5.9}  We have
$$\l\{w^3+\l\lfloor\f{x^3}2\r\rfloor+\l\lfloor\f{y^3}3\r\rfloor+\l\lfloor\f{z^3}4\r\rfloor:\ x,y,z\in\N\r\}=\N$$
and
$$\l\{w^3+\l\lfloor\f{x^3}2\r\rfloor+\l\lfloor\f{y^3}4\r\rfloor+\l\lfloor\f{z^3}8\r\rfloor:\ x,y,z\in\N\r\}=\N.$$
\endproclaim

Our following conjecture is a natural extension of Goldbach's Conjecture.

\proclaim{Conjecture 5.10} For any positive integers $a$ and $b$ with $a+b>2$, any integer $n>2$ can be written as $\lfloor p/a\rfloor+\lfloor q/b\rfloor$
with $p$ and $q$ both prime.
\endproclaim
\Remark\ 5.6. In the case $\{a,b\}=\{1,2\}$, Conjecture 5.10 reduces to Lemoine's Conjecture which states that any odd number greater than 5 can be written as $p+2q$
with $p$ and $q$ both prime. In the case $a=b=2$, Conjecture 5.10 reduces to the Goldbach Conjecture.
\medskip

Let us conclude this paper with one more conjecture.
\proclaim{Conjecture 5.11} Let
$$\align S=&\l\{\l\lfloor\f{x}9\r\rfloor:\ x-1\ \t{and}\ x+1\ \t{are twin prime}\r\}
\\=&\l\{\l\lfloor\f{x}3\r\rfloor:\ 3x-1\ \t{and}\ 3x+1\ \t{are twin prime}\r\}.
\endalign$$
Then, any positive integer can be written as the sum of two distinct elements of $S$ one of which is even.
Also, any positive integer can be expressed as the sum of an element of $S$ and a positive generalized pentagonal number.
\endproclaim
\Remark\ 5.7. Clearly either of the two assertions in Conjecture 5.11 implies the Twin Prime Conjecture.

\proclaim{Conjecture 5.12} Any integer $n>1$ can be written as $x^2+y^2+\varphi(z^2)$ with $x,y\in\N$, $z\in\Z^+$, and $\max\{x,y\}$ or $z$ prime.
Also, any $n\in\Z^+$ can be written as $x^3+y^2+T_z$ with $x,y\in\N$ and $z\in\Z^+$.
\endproclaim
\Remark\ 5.8. We have verified this for all $n=1,\ldots,10^5$. See [S15c, A262311 and A262813] for related data.

\proclaim{Conjecture 5.13} Any integer $m$ can be written as $x^4-y^3+z^2$ with $x,y,z\in\Z^+$.
\endproclaim
\Remark\ 5.9. We have verified this for all $m\in\Z$ with $|m|\ls 10^5$, see [S15c] for related data.
For example, 
$$0=4^4-8^3+16^2,\ \ 6= 36^4 - 139^3 + 1003^2,\ \ \t{and}\  11019 = 4325^4 - 71383^3 + 3719409^2.$$

\proclaim{Conjecture 5.14} Any $n\in\N$ can be written as $w^2+x^3+y^4+2z^4$ with $w,x,y,z\in\N$.
Also, any $n\in\N$ can be written as $w^2+2x^2+y^3+2z^3$ with $w,x,y,z\in\N.$
\endproclaim
\Remark\ 5.10. We have verified this for all $n=1,\ldots,4\times10^6$, see [S15c, A262827 and A262857] for related data.

\enddocument